\documentclass[11pt]{article}
\usepackage{amssymb,amsfonts,amsmath,amsthm}
\usepackage{epsfig}
\newtheorem{thm}{Theorem}[section]

\makeatletter \@addtoreset{equation}{section}

\def\qed{\hfill \rule{4pt}{7pt}}
\def\pf{\noindent {\it Proof.}\quad}
\title{\bf\LARGE Extended Zeilberger's Algorithm for Identities on Bernoulli and Euler Polynomials}
\author{William Y. C. Chen$^{1}$ and
Lisa H. Sun$^{2}$
\date{Center for Combinatorics, LPMC-TJKLC\\
 Nankai
University, Tianjin 300071, P.R. China\\
\vskip 0.2 cm Email: $^1$chen@nankai.edu.cn,
$^2$sun@cfc.nankai.edu.cn}}

\begin{document}
\maketitle

{\noindent \bf Abstract.} We present a computer algebra approach to
proving identities on Bernoulli polynomials and Euler polynomials by
using the extended Zeilberger's algorithm given by Chen, Hou and Mu.
The key idea is to use the contour integral definitions of the
Bernoulli and Euler numbers to establish recurrence relations on the
integrands. Such recurrence relations have certain parameter free
properties which lead to the required identities without computing
the integrals.



\section{Introduction}

Bernoulli polynomials and Euler polynomials play fundamental roles
in various branches of mathematics including combinatorics, number
theory, special functions and analysis, see for example
\cite{dilcher96,miki78,sunpan06}. At first glance, the Bernoulli
numbers, Euler numbers, and the corresponding polynomials do not
seem to fall in the framework of hypergeometric identities. The
powerful algorithm of Zeilberger \cite{zeil91} does not look like
the right mechanism to handle the Bernoulli and Euler numbers or
polynomials.

However, as  will be seen, the Cauchy contour integral
representations of the Bernoulli numbers and Euler numbers make it
possible to transform identities on these numbers and polynomials
into identities on hypergeometric sums. In order to avoid the
computation of the contour integrals, it is desirable to derive
recurrence relations of the hypergeometric summands with certain
parameter free properties. At the first trial, one finds it quite
disappointing that the recurrence relations given by Zeilberger's
algorithm seldom have the desired parameter free properties.
Nevertheless, this drawback can be overcome by using an extended
version of Zeilberger's algorithm. Paule \cite{paule05} first
noticed that Zeilberger's algorithm can be extended to derive mixed
recurrence relations for a hypergeometric term $F(n,m_1,m_2,\ldots,
m_r,k)$, where $r\geq 1$ and $m_1, m_2, \ldots, m_r$ are parameters.
Chen, Hou and Mu \cite{chenhm08} further extended Paule's algorithm
to the case with additional parameter free properties. In this
paradigm, many identities on Bernoulli and Euler numbers and
polynomials can be verified. Moreover, new identities can be
discovered from the recurrences generated by the original
Zeilberger's algorithm without the consideration of the parameter
free properties.

\section{Background}

Let us recall the background on Bernoulli and Euler numbers and
polynomials.  Let $\mathbb{N}=\{0,1,2,\cdots\}$ and
$\mathbb{Z}^{+}=\{1,2,3\cdots\}$. The well-known Bernoulli numbers
and Euler numbers are defined by the generating functions
$$
\sum_{n=0}^\infty B_n \frac{z^n}{n!}=\frac{z}{e^z-1} \quad
\mbox{and}\quad \sum_{n=0}^\infty E_n
\frac{z^n}{n!}=\frac{2e^z}{e^{2z}+1}.
$$
By the Cauchy integral formula, we have the contour integral
definitions of the Bernoulli numbers and the Euler numbers
\begin{align}
B_n&=\frac{n!}{2\pi i}\oint \frac{z}{e^z-1}\frac{d
z}{z^{n+1}},\label{condefi}\\[4pt]
E_n&=\frac{n!}{2\pi i}\oint \frac{2e^z}{e^{2z}+1}\frac{d
z}{z^{n+1}},\label{condefieuler}
\end{align}
where the contour encloses the origin, has radius less than $2\pi$
(to avoid the poles at $\pm 2\pi i$), and is traversed in a
counterclockwise direction. Actually, as will be seen, there will be
no need to compute the contour integrals, and one can formally treat
the contour integrals as linear operators. The integral
representation plays a crucial role in connecting the Bernoulli
numbers and Euler numbers to hypergeometric terms.

The Bernoulli numbers are also given by the following recursion
\begin{equation}\label{recber}
\sum_{k=0}^n {n+1\choose k}B_k=0, \quad n>0,
\end{equation}
with $B_0=1$. The Bernoulli numbers are rational and it is well
known that $B_{2n+1}=0$ for $n\geq 1$. The first few values of the
Bernoulli numbers are as follows
$$
B_0=1,\ B_1=-\frac{1}{2},\ B_2=\frac{1}{6},\ B_4=-\frac{1}{30}, \
B_6= \frac{1}{42}.
$$
For the Euler numbers, $E_{2n+1}=0$ for $n\geq 0$.

The Bernoulli polynomials and Euler polynomials can be defined by
the generating functions
$$
\sum_{n=0}^\infty B_n(x)\frac{t^n}{n!}=\frac{te^{xt}}{e^t-1} \quad
\mbox{and} \quad \sum_{n=0}^\infty E_n(x)
\frac{t^n}{n!}=\frac{2e^{xt}}{e^t+1}.
$$
Clearly  $B_n=B_n(0)$ and $E_n=2^nE_n(\frac{1}{2})$. The polynomials
$B_n(x)$ and $E_n(x)$  obey the following relations
\begin{align}
B_n(x)&=\sum_{k=0}^n {n\choose k}x^{n-k}B_k \label{berpol},\\
E_n(x)&=\sum_{k=0}^n {n\choose k} \Big(x-\frac{1}{2}\Big)^{n-k}
\frac{E_k}{2^k}.\label{eulerpol}
\end{align}

We will need the following basic properties of $B_n(x)$ and
$E_n(x)$. Lehmer \cite{leh88} showed that the Bernoulli polynomials
satisfy the relations $B_n(1)=(-1)^nB_n(0)$ and
\begin{equation}\label{berminus}
B_n(1-x)=(-1)^n B_n(x).
\end{equation}
Similarly, the Euler polynomials satisfy the relation
\begin{equation}\label{eulerminus}
E_n(1-x)=(-1)^n E_n(x).
\end{equation}
It is well known that the Bernoulli and Euler polynomials have the
following binomial expansions
\begin{equation}\label{bereulexpan}
B_n(x+y)=\sum_{k=0}^n {n \choose k} B_k(x) y^{n-k} \quad \mbox{and}
\quad E_n(x+y)=\sum_{k=0}^n {n \choose k} E_k(x) y^{n-k}.
\end{equation}

These basic properties will be needed for the computation of initial
values for the recurrence relations derived by our algorithm.

\section{The Algorithm}

In this section, we present an approach to proving Bernoulli number
identities by using the extended Zeilberger's algorithm, and we will
use an example to describe the four steps of our algorithm. The
original Zeilberger's algorithm is devised to find recurrence
relations of the summation $\sum_k F(n,k)$ by solving the equation
\[
a_0(n) F(n,k) +a_1(n)F(n+1, k) + \cdots + a_J(n)F(n+J,k)= G(n,k+1)-G(n,k),
\]
where $F(n,k)$ is a hypergeometric term in $n$ and $k$, $a_i(n)$ are
polynomials in $n$ and are $k$-free, $G(n,k)/F(n,k)$ is a rational
function in $n$ and $k$.  It is known that Zeilberger's algorithm
can be applied to summands with parameters in order to establish
multiple index recurrence relations, for example, see \cite[Section
4.3.1]{andpausch05} and \cite{paule05}. Recently, Chen, Hou and Mu
\cite{chenhm08} have found an extension of Zeilberger's algorithm to
summations of hypergeometric terms $\sum_k F(n,m_1, m_2, \ldots,
m_r,k)$, where $r \geq 1$ and $m_1,m_2,\ldots, m_r$ are parameters.
In fact, there are cases when the extended Zeilberger's algorithm
becomes more efficient than the original form, see \cite{chenhm08}.
 We will not give a rigorous description of the extended Zeilberger's algorithm,
since it will become apparent when it is being used.

For example, let us consider an identity of Gessel \cite[Lemma
7.2]{gessel03}.
\begin{thm} We have
\begin{equation}\label{gesselber}
\sum_{k=0}^m {m\choose k} B_{n+k}=(-1)^{m+n}\sum_{k=0}^n {n\choose
k} B_{m+k},
\end{equation}
where $m$ and $n$ are nonnegative integers.
\end{thm}

{\noindent \bf Proof.} To justify the above identity, we aim to find
recurrence relations for both sides. If they agree with each other,
then the equality is established by considering the initial values.
There are three steps to compute the recurrence relations for the
above summations. We will give detailed steps for the left hand side
of \eqref{gesselber}.

{\noindent \bf Step 1.} Extract the hypergeometric sum from the
Cauchy integral formula.

      Denote the left  hand side of \eqref{gesselber} by $L(n,m)$.
       By the contour integral formula for $B_n$, we have
\begin{align*}
L(n,m)&=\frac{1}{2\pi i} \oint \frac{1}{e^z-1} \sum_{k=0}^m
{m\choose k} \frac{(n+k)!}{z^{n+k}}\ dz.
\end{align*}
Denote the summand in the above integral by
$$
C(n,m,k)={m\choose k} \frac{(n+k)!}{z^{n+k}},
$$
and let \[ S(n,m)=\sum\limits_{k=0}^m C(n,m,k).\]

{\noindent \bf Step 2.} Construct an extended
telescoping equation with a shift  on the parameter $m$
of the summand $C(n,m,k)$, and solve this equation by
the extended Zeilberger's algorithm.

 Set the hypergeometric term
\begin{equation}\label{object}
F(n,m,k)=b_0C(n,m,k)+b_1C(n,m+1,k)+b_2C(n+1,m,k)+b_3 C(n+1,m+1,k),
\end{equation}
where $b_i$'s are $k$-free rational functions of $n$ and $m$,
namely, $k$ does not appear in $b_i$'s. Moreover, we require that
the rational functions $b_i$'s are independent of the variable $z$.

By Gosper's algorithm, it is easy to check that
$C(n,m,k)=z_{k+1}-z_{k}$ has no hypergeometric solution for $z_k$.
Moreover, since the Bernoulli numbers are not $P$-recursive, one
sees that  Zeilberger's algorithm does not work in this case.
Instead, we will try to solve the equation
\begin{equation}\label{recgosper}
F(n,m,k)=G(n,m,k+1)-G(n,m,k)
\end{equation}
where  $G(n,m,k)$ is a hypergeometric term. By Gosper's algorithm,
we get
\begin{equation}\label{normalfac}
r(k)=\frac{F(n,m,k+1)}{F(n,m,k)}=\frac{a(k)}{b(k)}\frac{c(k+1)}{c(k)},
\end{equation}
where
\begin{align*}
a(k)&=(m-k+1)(n+k+1),\\
b(k)&=z(k+1),\\
c(k)&=b_2k^2+(b_2n-b_2m+b_0z-b_3m-b_3)k-b_0z(m+1)\\
&\quad -b_1z(m+1)-b_2(n+1)(m+1)-b_3(n+1)(m+1).
\end{align*}

Assume that $G(n,m,k)=y(k)F(n,m,k)$, where $y(k)$ is an unknown
rational function of $k$. Substituting $y(k)F(n,m,k)$ for $G(n,m,k)$
in \eqref{recgosper} reveals that $y(k)$ satisfies
\begin{equation*}
r(k)y(k+1)-y(k)=1.
\end{equation*}
Substituting the factorization \eqref{normalfac} into the above
equation,  and setting \[ x(k)={y(k)c(k)\over b(k-1)},\] then
Zeilberger's algorithm reduces the problem further to that of
finding polynomial solutions (see \cite[Theorem 5.2.1]{aeqb}) of the
following equation
\begin{equation}\label{equgosper}
a(k)x(k+1)-b(k-1)x(k)=c(k).
\end{equation}
Notice that the coefficients $a(k)$ and $b(k)$ are independent of
the unknowns $b_i$'s, and $c(k)$ is a linear combination of $b_i$'s.
One can estimate the degree of the polynomial $x(k)$, as in Gosper's
algorithm. In this case, $x(k)$ is of degree $0$. Assume that
$x(k)=a_0$. Then the equation \eqref{equgosper} becomes
\begin{align*}
&(-a_0-b_2)k^2+\big(mb_2-nb_2+mb_3-b_0z-(n-m)a_0-a_0z+b_3\big)k\\
&+\big((n+1)(m+1)a_0+b_0z+nb_3+b_3nm+b_0zm+b_2nm+nb_2\\
&+b_1zm+b_1z+b_2+mb_2+mb_3+b_3\big)=0.
\end{align*}
By setting the coefficient of each power of $k$ to zero, we get a
system of linear equations in  $a_0$ and $b_i$'s. Note that in the
solution of this system, $a_0$ and $b_i$'s may contain the variable
$z$. To prevent $z$ from appearing in $b_i$'s, we should go one step
further to impose that the coefficient of any positive power of $z$
in $b_i$'s is zero. This may also lead to additional equations.
Combining all these equations, if we can find a nonzero solution,
then take this solution to the next step. Otherwise, we may try
recurrences of higher order. In this case, we get a nonzero solution
$a_0=-1, b_0=1, b_1=-1, b_2=1, b_3=0$. Note that in general the
$b_i$'s are polynomials in $n$ and $m$.

{\noindent \bf Step 3.} Compute the recurrence for $L(n,m)$.

By Step 2, the solution of $a_0, b_0, \ldots, b_3$ leads to the
following telescoping equation
\begin{equation}
C(n,m,k)-C(n,m+1,k)+C(n+1,m,k)=G(n,m,k+1)-G(n,m,k),
\end{equation}
where
\begin{equation}
G(n,m,k)=\frac{m!(n+k)!}{(k-1)!(m-k+1)!z^{n+k}}.
\end{equation}
Summing the above recurrence over $k$ from $0$ to $m+1$, we obtain
$$
S(n,m)-S(n,m+1)+S(n+1,m)=0.
$$
Substituting the above recurrence relation to the contour integral
definition of $B_n$, we  find that $L(n,m)$ satisfies
$$
L(n,m)-L(n,m+1)+L(n+1,m)=0.
$$

By the same procedure, we see that the right hand side of
\eqref{gesselber}, denoted by $R(n,m)$, satisfies the same
recurrence relation as $L(n,m)$, namely,
$$
R(n,m)-R(n,m+1)+R(n+1,m)=0.
$$

\noindent {\bf Step 4.} Verify initial values.

 By considering the parity of $B_m$, we see
that $(-1)^m B_m=B_m$ unless $m=1$. Therefore $L(0,1)=R(0,1)=1/2$
and
 $L(0,m)=R(0,m)=B_m$ for $m\neq 1$. This completes the proof. \qed

It is known that the Bernoulli numbers and Euler numbers are not
$P$-recursive, see \cite{fgs05}. Roughly speaking, this fact implies
that the original Zeilberger's algorithm is not applicable to derive
a recurrence relation of any order for summations involving
Bernoulli numbers.  For this reason, the extended Zeilberger's
algorithm becomes necessary, and it also suggests that in the study
of $P$-recursiveness of a polynomial sequence with parameters it is
likely that one can get a recurrence relation with polynomial
coefficients even for the sequence is not $P$-recursive, as long as
one allows shifts on the parameters.

\section{Bernoulli Number Identities}

In this section, we give several examples of proving identities on
Bernoulli numbers by using the extended Zeilberger's algorithm.

The first example is the extension of Kaneko's identity given by
Momiyama \cite{momiyama01}. It was proved by using a $p$-adic
integral over $\mathbb{Z}_p$. The Kaneko identity is stated as
follows
 \cite{kaneko95}
\begin{equation}\label{kanide}
\sum_{k=0}^{n+1} {n+1\choose k}\tilde{B}_{n+k}=0,
\end{equation}
where $\tilde{B}_{n}=(n+1)B_n$.

While our approach does not directly apply to Kaneko's identity
because it has no parameters, we can deal with Momiyama's identity
which reduces Kaneko's identity by setting $m=n$.

\begin{thm}[Momiyama's identity]
\begin{align}\label{momiyama}
(-1)^m\sum_{k=0}^m {m+1 \choose k} (n+k+1) B_{n+k} =-(-1)^n
\sum_{k=0}^n {n+1\choose k} (m+k+1) B_{m+k},
\end{align}
where $m$ and $n$ are integers and $m+n>0$.
\end{thm}

\pf Denote the left hand side and the right hand side of
\eqref{momiyama} by $L(n,m)$ and $R(n,m)$, respectively. By the
contour integral definition of the Bernoulli numbers, we have
\begin{equation*}
L(n,m) =\frac{1}{2\pi i} \oint \frac{1}{e^z-1} \Bigg(\sum_{k=0}^m
(-1)^m{m+1 \choose k} (n+k+1)\frac{(n+k)!}{z^{n+k}} \Bigg) dz.
\end{equation*}
Denote the summand in the above summation  by $F(n,m,k)$, that is,
\[
F(n,m,k)=(-1)^m{m+1 \choose k} (n+k+1)\frac{(n+k)!}{z^{n+k}}.
\]
Applying the extended Zeilberger's algorithm to $F(n,m,k)$ and
assuming that the output is independent of $z$, we obtain
\begin{equation} \label{fnmk}
F(n,m,k)+F(n,m+1,k)+F(n+1,m,k)=G(n,m,k+1)-G(n,m,k),
\end{equation}
where
$$
G(n,m,k)=\frac{(-1)^m(m+1)!(n+k+1)!}{(k-1)!(m+2-k)!z^{n+k}}.
$$
Summing the telescoping equation \eqref{fnmk} over $k$ from $0$ to
$m$,  we are led to the following recurrence relation for $L(n,m)$
$$
L(n,m)+L(n,m+1)+L(n+1,m)=-(-1)^m(n+m+2)B_{n+m+1}.
$$
Similarly, we find that $R(n,m)$ also satisfies
$$
R(n,m)+R(n,m+1)+R(n+1,m)=(-1)^n (n+m+2) B_{n+m+1}.
$$
Considering the parity of $B_n$, it is easy to see that
$$
\big((-1)^m+(-1)^n\big)(n+m+2)B_{n+m+1}=0.
$$
Therefore, both sides of Momiyama's identity \eqref{momiyama}
 satisfy the same recurrence relation.

To compute the initial values, setting $m=0$ we get
$L(n,0)=(n+1)B_n$. It follows from the recursion \eqref{recber}
that
$$
\sum_{k=0}^{n} {n+1 \choose k} B_k=\sum_{k=0}^{n} {n\choose
k}B_k+\sum_{k=0}^n {n \choose k-1}B_k=0.
$$
On the other hand, for $n\not=1$, we have
\begin{align*}
R(n,0)&=-(-1)^n \sum_{k=0}^n {n+1 \choose k} (k+1) B_k\\
&=-(-1)^n \Bigg(\sum_{k=0}^n{n+1 \choose k} k B_k+\sum_{k=0}^n {n+1
\choose k} B_k\Bigg)\\
&=-(-1)^n (n+1) \sum_{k=0}^n{n \choose k-1} B_k \\
&=(-1)^n (n+1)B_n=(n+1)B_n.
\end{align*}
It is easily checked that $L(1,0)=R(1,0)=-1$. So  we  deduce that
$L(n,0)=R(n,0)$ for all $n\geq 0$. This  completes the proof. \qed

The following identity is due to Gessel and Viennot \cite{gesvie}.

\begin{thm}[Gessel-Viennot]
\begin{equation}\label{beride}
\sum_{j=0}^{\lfloor(k-1)/2\rfloor} \frac{1}{k-j} {2k-2j \choose
k+1}{2n+1 \choose 2j+1} B_{2n-2j}=\frac{2n+1}{2k-2n+1} {2k-2n+1
\choose k+1}, \ n<k.
\end{equation}
\end{thm}

\pf Denote the left hand side and the right hand side of the above
identity by $L(n,k)$ and $R(n,k)$, respectively. So we get
\begin{align*}
L(n,k) &=\frac{1}{2\pi i} \oint \frac{1}{e^z-1}
\Bigg(\sum_{j=0}^{\lfloor(k-1)/2\rfloor} \frac{1}{k-j} {2k-2j
\choose k+1}{2n+1 \choose 2j+1}\frac{(2n-2j)!}{z^{2n-2j}}\Bigg) dz.
\end{align*}
Let $$ F(n,k,j)=  \frac{1}{k-j} {2k-2j \choose k+1}{2n+1 \choose
2j+1}\frac{(2n-2j)!}{z^{2n-2j}}.
$$
Applying the extended Zeilberger's algorithm, we get the following
recurrence
\begin{align*}
2(n+1)&(2n+3)F(n,k,j)+2(k+2)(2k+3)F(n+1,k+1,j)\\[3pt]
&-(k+2)(k+3)F(n+1,k+2,j)=G(n,k,j+1)-G(n,k,j),
\end{align*}
where
$$
G(n,k,j)=\frac{4j(2n+3)!(2k-2j+1)!}{(k-2j+1)!(2j)!(k+1)!z^{2n-2j+2}}.
$$
By summing the above telescoping equation over $j$, we obtain the
following recurrence relation for $L(n,k)$
\begin{equation}\label{rec}
2(n+1)(2n+3)L(n,k)+2(k+2)(2k+3)L(n+1,k+1)-(k+2)(k+3)L(n+1,k+2)=0.
\end{equation}
It is easy to check that  $R(n,k)$ also satisfies the above
recurrence relation. Since $n<k$, we can define $L(n,n)=R(n,n)=0$
for $n \neq 0$. It is also easy to verify the initial conditions
\[ L(0,k)=R(0,k)=\frac{1}{2k+1}{2k+1 \choose k+1}.\] This completes the
proof. \qed

It should be noted that the recurrence relation \eqref{rec} for
$L(n,k)$ was derived by Jacobi \cite{jacobi34} in 1834, see Gessel
and Viennot \cite{gesvie}.

The next identity is due to Gelfand \cite{gelfand68}.

\begin{thm} We have
\begin{equation}\label{gelfandide}
(-1)^{n-1} \sum_{k=1}^m {m\choose k-1}
\frac{B_{n+k}}{n+k}+(-1)^{m-1} \sum_{k=1}^n {n\choose k-1}
\frac{B_{m+k}}{m+k}=\frac{m!n!}{(m+n+1)!},
\end{equation}
provided that the integers $m,n\geq 0$ are not both zero.
\end{thm}

\pf Denote the left and right hand sides of the above identity
\eqref{gelfandide} by $L(n,m)=S(n,m)+T(n,m)$ and $R(n,m)$,
respectively, where $S(n,m)$ and $T(n,m)$ are the first and second
sums of $L(n,m)$. Note that
\begin{align*}
S(n,m)&=\frac{1}{2\pi i} \oint \frac{1}{e^z-1} \bigg(\sum_{k=1}^m
(-1)^{n-1} {m\choose k-1} \frac{(n+k)!}{(n+k)z^{n+k}}\bigg) dz,\\[5pt]
T(n,m)&=\frac{1}{2\pi i} \oint \frac{1}{e^z-1} \bigg(\sum_{k=1}^n
 (-1)^{m-1} {n\choose k-1} \frac{(m+k)!}{(m+k)z^{m+k}}\bigg) dz.
\end{align*}
Denote the summand in $S(n,m)$ by $F(n,m,k)$, and by the extended
Zeilberger's algorithm, we obtain
$$
F(n,m,k)-F(n,m+1)-F(n+1,m)=G(n,m,k+1)-G(n,m,k),
$$
where
$$
G(n,m,k)=(-1)^{n-1}\frac{m!(n+k-1)!}{(k-2)!(m+2-k)!z^{n+k}}.
$$
Summing the above telescoping equation over $k$ from $1$ to $m$, we
get a recurrence for $S(n,m)$
$$
S(n,m)-S(n,m+1)-S(n+1,m)=(-1)^n \frac{B_{m+n+1}}{m+n+1}.
$$
By the same procedure, or  by the symmetric property
$T(n,m)=S(m,n)$, we find that
$$
T(n,m)-T(n,m+1)-T(n+1,m)=(-1)^m \frac{B_{m+n+1}}{m+n+1}.
$$
With the aid of the property  $B_{2n+1}=0$ for $n\ge 1$, we have
$$
L(n,m)-L(n,m+1)-L(n+1,m)=\big((-1)^m+(-1)^n\big)\frac{B_{m+n+1}}{m+n+1}=0.
$$
It is easy to verify that $R(n,m)$ also satisfies the above
recurrence relation. To check the initial values, we have
\begin{align*}
L(n,0)=0-\sum_{k=1}^n {n\choose k-1}\frac{B_k}{k}
=-\frac{1}{n+1}\sum_{k=1}^n {n+1 \choose k} B_k
=\frac{1}{n+1}=R(n,0).
\end{align*}
This completes the proof. \qed

Agoh and Dilcher \cite[Theorem 2.1]{agodil07} obtained a convolution
identity for Bernoulli numbers. By the extended Zeilberger's
algorithm and Woodcock's identity \eqref{woodide}, we can give a
direct proof of this result which is restated  in the following
equivalent form.

\begin{thm} Let $m,n,k \geq 0$ be integers, with $m$ and $k$ not
both zero. Then
\begin{align}\label{agodilide}
\sum_{j=0}^n {n\choose j} B_{k+j}
B_{m+n-j}&=-\frac{k!m!}{(m+k+1)!}\big(n+\delta(m,k)(m+k+1)\big)B_{m+n+k}\nonumber\\
&+\sum_{r=0}^{m+k} (-1)^r
\frac{B_{m+k+1-r}}{m+k+1-r}(-1)^k{k+1\choose
r}\Big(\frac{k+1-r}{k+1}n-\frac{rm}{k+1}\Big)B_{n+r-1}\nonumber\\
&+\sum_{r=0}^{m+k} (-1)^r
\frac{B_{m+k+1-r}}{m+k+1-r}(-1)^m{m+1\choose
r}\Big(\frac{m+1-r}{m+1}n-\frac{rk}{m+1}\Big)B_{n+r-1},
\end{align}
where $\delta(m,k)=0$ when $m=0$ or $k=0$, and $\delta(m,k)=1$
otherwise.
\end{thm}

\pf Let $L(n,m,k)$ and $R(n,m,k)$ denote the left hand side and the
right hand side of the above identity \eqref{agodilide},
respectively. Our approach leads to the  recurrence relation
\begin{equation}
S(n,m+1,k)-S(n+1,m,k)+S(n,m,k+1)=0,
\end{equation}
where $m\neq 0$ and $k\neq 0$. Considering the parity of the
Bernoulli numbers, we have $(-1)^kB_k=B_k$ for $k \neq 1$. The known
convolution identity on Bernoulli numbers
\begin{equation}
\sum_{k=0}^n {n\choose k}B_k B_{n-k}=-nB_{n-1}-(n-1)B_n,\quad n\geq
1
\end{equation} yields that
\begin{align*}
L(0,m,1)=&B_1B_m=-\frac{1}{2}B_m,\\[3pt]
R(0,m,1)=&-\frac{1}{m+1}B_{m+1}+\sum_{r=0}^{m+1} (-1)^r
\frac{B_{m+2-r}}{m+2-r}(-1){2\choose
r}\Big(-\frac{rm}{2}\Big)B_{r-1}\\
&+\sum_{r=0}^{m+1} (-1)^r \frac{B_{m+2-r}}{m+2-r}(-1)^m{m+1\choose
r}\Big(-\frac{r}{m+1}\Big)B_{r-1}\\
=&
-\frac{1}{m+1}B_{m+1}-\frac{m}{m+1}B_{m+1}+B_mB_1+\frac{(-1)^m}{m+1}
\sum_{r=0}^{m} (-1)^r {m+1\choose r} B_{m+1-r}B_r\\
=&-B_{m+1}+B_mB_1+\frac{(-1)^m}{m+1}\bigg( \sum_{r=0}^{m+1}
{m+1\choose r} B_{m+1-r}B_r-2(m+1)B_mB_1-B_{m+1}\bigg) \\
=&-B_{m+1}+B_mB_1-(-1)^mB_{m+1} \\
=&-\frac{1}{2}B_m.
\end{align*}
This gives the proof for \eqref{agodilide} when $m\neq 0$ and $k\neq
0$.

Moreover, if $m=0$ or $k=0$, we can simplify the identity to an
equivalent form of a known identity discovered by Woodcock
\cite{woodcock79}
\begin{equation}\label{woodide}
\frac{1}{m}\sum_{k=1}^m (-1)^k{m \choose k}
B_{m-k}B_{n-1+k}=\frac{1}{n} \sum_{k=1}^n (-1)^k{n \choose k}
B_{n-k}B_{m-1+k}.
\end{equation}
This completes the proof.  \qed

\section{Bernoulli Polynomial Identities}

In this section, we show that our approach is also valid for proving
identities on Bernoulli polynomials. We will explain how this method
works by considering an identity due to Sun \cite{sun031}.

\begin{thm} We have
\begin{equation}\label{ezberpol}
(-1)^k\sum_{j=0}^k {k\choose j} x^{k-j} B_{l+j}(y)
=(-1)^l\sum_{j=0}^l {l\choose j} x^{l-j} B_{k+j}(z),
\end{equation}
provided that $x+y+z=1$.
\end{thm}

\pf Denote both sides of the above equation by $L(k,l)$ and
$R(k,l)$, respectively.  We have
\begin{align*}
L(k,l)=\frac{1}{2\pi i} \oint \frac{1}{e^u-1}\Bigg(\sum_{j=0}^k
\sum_{h=0}^{l+j} (-1)^k{k\choose j} {l+j \choose h}x^{k-j} y^{l+j-h}
\frac{h!}{u^h}\Bigg) du.
\end{align*}
Let $F(k,l,h,j)$ denote the summand in the above integral, that is,
\[ F(k,l,h,j)=(-1)^k{k\choose j} {l+j \choose h}x^{k-j} y^{l+j-h}
\frac{h!}{u^h}.\] Applying the extended Zeilberger's algorithm to
$F(k,l,h,j)$ with the assumption that the output is independent of
the variables $u$ and $h$, we arrive at the relation
\begin{equation}\label{recezpol}
xF(k,l,h,j)+F(k+1,l,h,j)+F(k,l+1,h,j)=G(k,l,h,j+1)-G(k,l,h,j),
\end{equation}
where \[  G(k,l,h,j)={ xj\over k-j+1} F(k,l,h,j) .\]

Summing both sides of  \eqref{recezpol} over $h$ and $j$ gives the
recurrence relation
$$
xL(k,l)+L(k+1,l)+L(k,l+1)=0.
$$
Similarly, it can be shown that $R(k,l)$ satisfies the same
recurrence relation. It remains to check  the initial values
\begin{align*}
L(0,l)&=B_l(y),\\
R(0,l)&=(-1)^l \sum_{j=0}^l {l\choose j} x^{l-j} B_j(z)=(-1)^l \sum_{j=0}^l {l\choose j} x^{l-j} (B+z)^j\\
&=(-1)^l (B+x+z)^l=(-1)^l B_l(x+z)=(-1)^l B_l (1-y)=B_l(y),
\end{align*}
as desired.
 \qed

It is worth noting that the extended Zeilberger's algorithm is
indeed efficient in deriving recurrence relations for multiple sums.
The next identity is given by Wu, Sun and Pan \cite{wusunpan04}.

\begin{thm} We have
\begin{align}\label{wspberpol}
&(-1)^m \sum_{k=0}^m {m+1 \choose k}(n+k+1)B_{n+k}(x)\nonumber\\
&+(-1)^n \sum_{k=0}^n {n+1\choose k} (m+k+1)
B_{m+k}(-x)\nonumber\\[3pt]
& \qquad =\,(-1)^m(n+m+1)(n+m+2)x^{n+m}.
\end{align}
\end{thm}

\pf Denote the two sums on the left hand side of (\ref{wspberpol})
by $S(n,m)$ and $T(n,m)$ respectively. Let $L(n,m)=S(n,m)+T(n,m)$,
and let $R(n,m)$ denote the right hand side of (\ref{wspberpol}).
Write
\begin{align*}
S(n,m)=\frac{1}{2\pi i} \oint \frac{1}{e^z-1} \Bigg(
\sum_{k=0}^m\sum_{j=0}^{n+k} (-1)^m{m+1 \choose k}(n+k+1){n+k
\choose j} x^{n+k-j} \frac{j!}{z^j}\Bigg) dz.
\end{align*}
Denote the summand in the above expression by $F(n,m,k,j)$. Applying
the extended Zeilberger's algorithm with the assumption that the
output is independent of the parameters $z$ and $j$, we obtain that
\begin{align*}
F(n,m,k,j)+F(n+1,m,k,j)+F(n,m+1,k,j)=G(n,m,k+1,j)-G(n,m,k,j),
\end{align*}
where \[G(n,m,k,j)={k\over m-k+1} F(n,m,k,j).\] By summing the above
telescoping equation over $j$ from $0$ to $n+k$ and $k$ from $0$ to
$m+1$, we deduce that
\begin{equation}\label{rect1}
S(n,m)+S(n+1,m)+S(n,m+1)=(-1)^{m+1}(n+m+2)B_{n+m+1}(x).
\end{equation}
From the symmetry property it follows that $T(n,m)(x)=S(m,n)(-x)$.
This leads to the following recurrence relation for $T(n,m)$
\begin{equation}\label{rect2}
T(n,m)+T(n,m+1)+T(n+1,m)=(-1)^{n+1}(n+m+2)B_{n+m+1}(-x).
\end{equation}
Adding  \eqref{rect1} to \eqref{rect2}, we derive a recurrence
relation satisfied by $L(n,m)$ \allowdisplaybreaks
\begin{align*}
&L(n,m)+L(n+1,m)+L(n,m+1)\\[4pt]
&=(-1)^{m+1}(n+m+2)B_{n+m+1}(x)+(-1)^{n+1}(n+m+2)B_{n+m+1}(-x)\\[4pt]
&=(-1)^{m+1}(n+m+2)\sum_{k=0}^{n+m+1} {n+m+1\choose k}
x^{n+m+1-k}B_k
\bigg(1+(-1)^{k+1}\bigg)\\[4pt]
&=2(-1)^{m+1}(n+m+2)\sum_{{k=0} \atop {k,odd}}^{n+m+1} {n+m+1\choose
k}
x^{n+m+1-k}B_k \\[4pt]
&=2(-1)^{m+1}(n+m+2)(n+m+1)x^{n+m}B_1\\[4pt]
&=(-1)^{m}(n+m+1)(n+m+2)x^{n+m}.
\end{align*}
It is easy to see that $R(n,m)$ satisfies the same recurrence
relation as $L(n,m)$. Based on the well-known identity for Bernoulli
polynomials
\begin{equation*}
 nx^{n-1}=\sum_{k=1}^n {n
\choose k}B_{n-k}(x)=\sum_{k=0}^{n-1} {n\choose k}B_k(x),
\end{equation*}
it is straightforward to verify that
\begin{align*}
L(n,-1)&=0+(-1)^n\sum_{k=0}^n {n+1\choose k}kB_{k-1}(-x)\\[4pt]
&=(-1)^n(n+1) \sum_{k=0}^{n-1} {n\choose k}B_{k}(-x)\\[4pt]
&=(-1)^n(n+1)n(-x)^{n-1}=-n(n+1)x^{n-1}=R(n,-1).
\end{align*}
This completes the proof. \qed

Note that the above identity \eqref{wspberpol} reduces to Momiyama's
identity \eqref{momiyama} by setting $x=0$. We also note that
integrating the identity \eqref{wspberpol} over $x$ and using the
Bernoulli number identity \eqref{gesselber}, one can derive the
following identity of Wu, Sun and Pan \cite{wusunpan04}
\begin{equation}\label{wupansunide}
(-1)^m \sum_{i=0}^m {m \choose i } B_{n+i}(x)=(-1)^n\sum_{j=0}^n {n
\choose j}B_{m+j}(-x).
\end{equation}

The following identity is derived by Sun \cite{sun031}.

\begin{thm}We have
\begin{align}\label{sunide}
(-1)^k&\sum_{j=0}^k {k\choose j} x^{k-j} \frac{B_{l+j+1}(y)}{l+j+1}
+(-1)^l\sum_{j=0}^l {l\choose j} x^{l-j}
\frac{B_{k+j+1}(z)}{k+j+1}=\frac{(-x)^{k+l+1}}{(k+l+1){k+l \choose
k}},
\end{align}
provided that $x+y+z=1$.
\end{thm}

\pf Let $L(k,l)$ and $R(k,l)$ denote the left hand side and the
right hand side of (\ref{sunide}), respectively. It can be shown
that
$$
xL(k,l)+L(k+1,l)+L(k,l+1)=0.
$$
It can also be shown that $R(k,l)$ satisfies the same recurrence
relation. To check the  initial conditions, we have
\begin{align*}
L(0,l)&=\frac{B_{l+1}(y)}{l+1}+(-1)^l\sum_{j=0}^l {l\choose j}
x^{l-j} \frac{B_{j+1}(z)}{j+1}\\[3pt]
&=\frac{B_{l+1}(y)}{l+1}+\frac{(-1)^l}{l+1}\sum_{j=0}^l {l+1\choose
l-j+1} x^{j} B_{l-j+1}(z)\\[3pt]
&=\frac{B_{l+1}(y)}{l+1}+\frac{(-1)^l}{l+1}\sum_{j=0}^l {l+1\choose
j} x^{j} (B+z)^{l-j+1}\\[3pt]
&=\frac{B_{l+1}(y)}{l+1}+\frac{(-1)^l}{l+1}(B+x+z)^{l+1}-\frac{(-1)^l}{l+1}x^{l+1}\\[3pt]
&=\frac{B_{l+1}(y)}{l+1}+\frac{(-1)^l}{l+1}B_{l+1}(1-y)-\frac{(-1)^l}{l+1}x^{l+1}\\[3pt]
&=\frac{B_{l+1}(y)}{l+1}-\frac{1}{l+1}B_{l+1}(y)-\frac{(-1)^l}{l+1}x^{l+1}\quad (\mbox{by\ }\ \eqref{berminus})\\[3pt]
&=\frac{(-x)^{l+1}}{l+1}=R(0,l),
\end{align*}
as desired. \qed

We remark that the above identity \eqref{sunide} reduces to
\eqref{ezberpol} by viewing $z=1-x-y$ as a function of $y$ and by
taking partial derivative with respect to $y$. It  also specializes
to \eqref{wupansunide} when setting $y\rightarrow x$ and
$z=-y\rightarrow -x$. Moreover, differentiating both sides of
\eqref{sunide} with respect to $y$ twice, we obtain the following
identity derived by Sun \cite{sun031}, which can be verified by our
approach. The proof is omitted.

\begin{thm}Suppose that $x+y+z=1$, then
\begin{align}\label{sunide2}
&(-1)^k \sum_{j=0}^k {k+1 \choose j} x^{k-j+1} (l+j+1) B_{l+j}(y)\nonumber\\
&+(-1)^l \sum_{j=0}^l {l+1 \choose j} x^{l-j+1} (k+j+1) B_{k+j}
(z)\nonumber\\
&\qquad =\,(-1)^k (k+l+2) (B_{k+l+1}(x+y)-B_{k+l+1}(y)).
\end{align}
\end{thm}

In \cite[Theorem 1.1]{sunpan06}, Sun and Pan find a symmetric
relation between products of the  Bernoulli polynomials.

\begin{thm}Let $n\in \mathbb{Z}^{+}$ and $x+y+z=1$. If $r+s+t=n$, then
\begin{align}\label{bersymm}
&r\sum_{k=0}^n (-1)^k {s\choose k}{t \choose n-k}
B_{n-k}(x)B_k(y)\nonumber \\
&+s\sum_{k=0}^n (-1)^k {t\choose k}{r \choose n-k}
B_{n-k}(y)B_k(z)\nonumber \\
& +t \sum_{k=0}^n (-1)^k {r\choose k}{s \choose n-k}
B_{n-k}(z)B_k(x)=0.
\end{align}
\end{thm}

\pf  Denote the three sums on the left hand side of the above
identity by $S(n,r,s)$, $T(n,r,s)$, $R(n,r,s)$ respectively. Since
$n=r+s+t$, $S(n,r,s)$ can be expressed as
\[ \Big(\frac{1}{2\pi i}\Big)^2 \oint \frac{1}{e^u-1} \oint
\frac{1}{e^v-1} \Bigg(\sum_{k,j,h} (-1)^k{s\choose k}{n-r-s \choose
n-k}{n-k\choose j}{k \choose h} \frac{j!}{u^j}\frac{h!}{v^h} r
x^{n-k-j} y^{k-h} \Bigg) dudv.
\]
Our approach yields the following recurrence relation
$$
(s+1)S(n,r+1,s)+(r+1)S(n,r,s+1)+(n-r-s-1)S(n,r+1,s+1)=0.
$$
Similarly, it can be shown that $T(n,r,s)$ and $R(n,r,s)$ satisfy
the same recurrence relation. Since $r+s+t=n$, we obtain that
$$
S(n,0,s)+T(n,0,s)+R(n,0,s)=(-1)^n s {n-s \choose
n}B_n(z)+(n-s){s\choose n} B_n(z)=0,
$$
and
$$
S(n,r,0)+T(n,r,0)+R(n,r,0)=r {n-r \choose
n}B_n(x)+(n-r)(-1)^n{r\choose n} B_n(x)=0.
$$
It follows that $S(n,r,s)+T(n,r,s)+R(n,r,s)$ is identically zero.
This completes the proof for all integers $r, s$ and $t$. Then by taking the left hand side of  \eqref{bersymm} as a polynomial in $r, s, t$, it follows that \eqref{bersymm} is true for all  $r,s,t$ such that $r+s+t=n$.    \qed

\section{Euler Number and Polynomial Identities}

In this section, we show how to prove identities on Euler numbers
and polynomials by using our approach. As the first example, we
consider the following  identity due to Wu, Sun and Pan
\cite{wusunpan04}.

\begin{thm} We have
\begin{equation}\label{wspeuler}
(-1)^m\sum_{k=0}^m {m \choose k} \frac{E_{n+k}}{2^{n+k}}=(-1)^n
\sum_{j=0}^n {n \choose j} E_{m+j}\left(-\frac{1}{2}\right),
\end{equation}
where $m$ and $n$ are nonnegative integers.
\end{thm}

\pf Denote the left and right hand sides of \eqref{wspeuler} by
$L(n,m)$ and $R(n,m)$, respectively. By the contour integral
definition of the Euler numbers \eqref{condefieuler} and the
relation \eqref{eulerpol}, we have
\begin{align*}
L(n,m)&=\frac{1}{2\pi i} \oint \frac{2e^z}{e^{2z}+1}
\Bigg(\sum_{k=0}^m (-1)^m {m \choose k}
\frac{(n+k)!}{2^{n+k}z^{n+k+1}}\Bigg) dz,\\[5pt]
R(n,m)&=\frac{1}{2\pi i} \oint \frac{2e^z}{e^{2z}+1}
\Bigg(\sum_{j=0}^n \sum_{k=0}^{m+j} (-1)^n{n\choose j} {m+j\choose
k} (-1)^{m+j-k} \frac{k!}{2^k z^{k+1}}\Bigg) dz.
\end{align*}
Denote the summands in the above two integrands by
\begin{align*}
S(m,n,k)&=(-1)^m {m \choose k} \frac{(n+k)!}{2^{n+k}z^{n+k+1}},\\[5pt]
T(m,n,k,j)&=(-1)^n{n\choose j} {m+j\choose k} (-1)^{m+j-k}
\frac{k!}{2^k z^{k+1}}.
\end{align*}
Applying the extended Zeilberger's algorithm, we obtain
\begin{align*}
&S(n,m,k)+S(n,m+1,k)+S(n+1,m,k)=G(n,m,k+1)-G(n,m,k),\\[4pt]
&T(n,m,k,j)+T(n,m+1,k,j)+T(n+1,m,k,j)=H(n,m,k,j+1)-H(n,m,k,j),
\end{align*}
where
$$
G(n,m,k)=\frac{k}{m+1-k}S(n,m,k),\quad
H(n,m,k,j)=\frac{j}{n+1-j}T(n,m,k,j).
$$
Therefore, $L(n,m)$ and $R(n,m)$ satisfy the same recurrence
$$
L(n,m)+L(n,m+1)+L(n+1,m)=0.
$$
Consequently, the identity \eqref{wspeuler} can be verified by
computing the initial values
\begin{align*}
L(0,m)&=(-1)^m\sum_{k=0}^m {m\choose k}\frac{E_k}{2^k} =\sum_{k=0}^m
{m\choose k} (-1)^{m-k} \frac{E_k}{2^k}
=E_m\Big(-\frac{1}{2}\Big)=R(0,m),
\end{align*}
as desired. \qed

Wu, Sun and Pan \cite{wusunpan04} also derived an identity by
substituting the Bernoulli polynomials in \eqref{wupansunide} with
 Euler polynomials. This identity can be verified by our approach.
 The proof is omitted.

\begin{thm} We have
\begin{align}\label{wspeulpol}
&(-1)^m \sum_{k=0}^m {m \choose k}E_{n+k}(x)=(-1)^n \sum_{k=0}^n
{n\choose k} E_{m+k}(-x).
\end{align}
\end{thm}

Note that differentiating both sides of the identity
\eqref{wspeulpol} with respect to $x$ leads to the following
identity also due to Wu, Sun and Pan \cite{wusunpan04}:
\begin{align}\label{wspeulpol2}
&(-1)^m \sum_{k=0}^m {m+1 \choose k}(n+k+1)E_{n+k}(x)\nonumber\\[4pt]
&+(-1)^n \sum_{k=0}^n {n+1\choose k} (m+k+1)
E_{m+k}(-x)\nonumber\\[4pt]
&=(-1)^m2(n+m+2)(x^{n+m+1}-E_{n+m+1}(x)).
\end{align}

The following identity is derived by Sun \cite{sun031}.

\begin{thm}
\begin{align}\label{suneuleride}
(-1)^k&\sum_{j=0}^k {k\choose j} x^{k-j} \frac{E_{l+j+1}(y)}{l+j+1}
+(-1)^l\sum_{j=0}^l {l\choose j} x^{l-j}
\frac{E_{k+j+1}(z)}{k+j+1}=\frac{(-x)^{k+l+1}}{(k+l+1){k+l \choose
k}},
\end{align}
provided that $x+y+z=1$.
\end{thm}

\pf  Denote the two sums in the left hand side of the above identity
by $S(k,l)$ and $T(k,l)$ respectively. Let $L(k,l)=S(k,l) + T(k,l)$,
and let $R(k,l)$ denote the right hand side of \eqref{suneuleride}.
By computation, we find
$$
x S(k,l)+S(k+1,l)+S(k,l+1)=0.
$$
Since $T(k,l)=S(l,k)$,
$$
xT(k,l)+T(k+1,l)+T(k,l+1)=0.
$$
Therefore,
$$
xL(k,l)+L(k+1,l)+L(k,l+1)=0.
$$
It is easy to check  that $R(k,l)$  satisfies the same recurrence
relation. To check the initial values, we have
\begin{align*}
L(0,l)&=\frac{E_{l+1}}{l+1}+(-1)^l \frac{1}{l+1} \sum_{j=0}^l {l+1
\choose j+1} x^{l-j} E_{j+1}(z)\\
&=\frac{E_{l+1}}{l+1}+(-1)^l \frac{1}{l+1} \sum_{j=1}^{l+1} {l+1
\choose j} x^{l+1-j} E_{j}(z)\\
&=\frac{E_{l+1}}{l+1}+(-1)^l \frac{1}{l+1}
\Big(E_{l+1}(x+z)-x^{l+1}\Big)\\
&=\frac{(-x)^{l+1}}{l+1}=R(0,l),
\end{align*}
as desired. \qed

Our approach  can also be applied to identities involving products
of the Euler polynomials and the Bernoulli polynomials. We take the
following identity of Sun and Pan \cite[Theorem 1.1]{sunpan06} as an
example. Note that in Sun and Pan's identity, the variables $r, s$ and $t$ should be real numbers.

\begin{thm}
Let $n\in \mathbb{Z}^{+}$, $r+s+t=n-1$ and $x+y+z=1$,
then
\begin{align}\label{bereuleride}
&\sum_{k=0}^n (-1)^k {r \choose k}{s \choose n-k} B_k(x)
E_{n-k}(z)\nonumber \\ &-(-1)^n \sum_{k=0}^n (-1)^k {r\choose k}{t
\choose
n-k}B_k(y)E_{n-k}(z)\nonumber \\
&\qquad =\frac{r}{2} \sum_{l=0}^{n-1} (-1)^l {s \choose l} {t\choose
n-1-l} E_l(y)E_{n-1-l}(x).
\end{align}
\end{thm}

\pf Denote the two sums on the left hand side of the above identity
\eqref{bereuleride} by  $S(n,r,s)$ and $T(n,r,s)$ respectively. Let
$L(n,r,s)=S(n,r,s)-T(n,r,s)$, and denote the right hand side of
\eqref{bereuleride} by $R(n,r,s)$. Note that
\begin{align*}
S(n,r,s)=\Bigg(\frac{1}{2\pi i}\Bigg)^2 \oint \frac{1}{e^u-1}\oint &
\frac{2e^v}{e^{2v}+1}   \Bigg(\sum_{k=0}^n \sum_{j=0}^k
\sum_{h=0}^{n-k}(-1)^k {r \choose k}{s \choose n-k} {k\choose
j}\\
&\times x^{k-j} \frac{j!}{u^j}{n-k\choose
h}\Big(z-\frac{1}{2}\Big)^{n-k-h}\frac{h!}{2^h v^{h+1}} \Bigg) dudv.
\end{align*}
Applying the  extended Zeilberger's algorithm, we have
$$
(s+1)S(n,r+1,s)+(r+1)S(n,r,s+1)+(n-s-r-2)S(n,r+1,s+1)=0.
$$
It can also be shown  that $T(n,r,s)$ and $R(n,r,s)$ satisfy the
same recurrence relation. Setting $r=0$, since $r+s+t=n-1$, it follows that
\begin{align*}
&L(n,0,s)=\bigg({s \choose n}-(-1)^n{n-1-s \choose n}\bigg)E_{n}(z)=0=R(n,0,s).
\end{align*}
To show that $L(n,r,0)=R(n,r,0)$, it is equivalent to verify
\begin{align}\label{bereuler2}
\sum_{k=0}^n (-1)^{n+k} {r
\choose k} {n-1-r \choose n-k}B_k(y)E_{n-k}(z)=(-1)^n{r \choose n}B_n(x)-\frac{r}{2} {n-1-r \choose n-1} E_{n-1}(x).
\end{align}
It is easy to see that both sides of the above identity satisfy the following recurrence relation
\[
(r+1)S(n, r)+(n-1-r)S(n,r+1)=0.
\]
 Then \eqref{bereuler2} can be proved by checking the initial case $r=0$.
This completes the proof of \eqref{bereuleride} for all integers $r, s$ and $t$.
Considering  both sides of  \eqref{bereuleride} as polynomials in $r, s, t$, we deduce that \eqref{bereuleride} holds for all  $r,s,t$ such that $r+s+t=n-1$. \qed

\section{Deriving Identities from Kaneko's Identity}

Applying the original Zeilberger's algorithm to a Bernoulli number
summation, we may obtain  a recurrence relation for the summand
which contains the integral variable $z$. Although such a recurrence
cannot be used to prove the Bernoulli number identity itself, it may
lead to another identity. For example, let us consider Kaneko's
identity \eqref{kanide}
\begin{equation*}
\sum_{k=0}^{n+1} {n+1\choose k}\tilde{B}_{n+k}=0,
\end{equation*}
where $\tilde{B}_{n}=(n+1)B_n$. From the recurrence obtained by
Zeilberger's algorithm, we can get the following generalization of
this identity.

\begin{thm} We have
\begin{equation}\label{kanger}
\sum_{k=0}^{n+3} {n+3\choose k} (n+k+3)(n+k+2)\tilde{B}_{n+k}=0.
\end{equation}
\end{thm}
 \pf Denote the left hand side of Kaneko's identity by $L(n)$. By the contour integral
definition of the Bernoulli numbers, we have \allowdisplaybreaks
\begin{align*}
L(n)
&=\sum_{k=0}^{n+1} {n+1\choose k} (n+k+1) B_{n+k}\\[3pt]
&=\frac{1}{2\pi i} \oint \frac{1}{e^z-1} \Bigg(\sum_{k=0}^{n+1}
{n+1\choose k} (n+k+1)\frac{(n+k)!}{z^{n+k}}\Bigg) dz.
\end{align*}
Denote the summation in the above integral by $S(n)$. Obviously,
\[L(n)=\frac{1}{2\pi i} \oint \frac{1}{e^z-1} S(n) dz=0\] for
all $n\geq 0$. Applying Zeilberger's algorithm, we get
\begin{align*}
z^2S(n+2)=2(n+3)(2n+5)S(n+1)+(n+2)(n+3)S(n).
\end{align*}
By integrating over $z$ on both sides of the above recurrence, it
follows that
\begin{align*}
& \frac{1}{2\pi i} \oint \frac{1}{e^z-1}  z^2S(n+2) dz\\
&= \frac{1}{2\pi i} \oint \frac{1}{e^z-1} \Bigg(\sum_{k=0}^{n+3}
{n+3\choose k} (n+k+3)\frac{(n+k+2)!}{z^{n+k}}\Bigg) dz\\
&=\sum_{k=0}^{n+3} {n+3\choose k} (n+k+3)(n+k+2)(n+k+1) B_{n+k}\\
&=\sum_{k=0}^{n+3} {n+3\choose k} (n+k+3)(n+k+2)\tilde{B}_{n+k}\\
&=2(n+3)(2n+5) \frac{1}{2\pi i} \oint \frac{1}{e^z-1}
 S(n+1) dz+(n+2)(n+3) \frac{1}{2\pi i} \oint \frac{1}{e^z-1}  S(n)
dz =0.
\end{align*}
This completes the proof. \qed

It should be mentioned that the above identity \eqref{kanger} is the special case $s=r=3$ of an identity of K.-W. Chen \cite{KWChen}:
\begin{equation} \label{kkk}
\sum_{k=0}^{n+r}{n+r\choose k}{n+r+k\choose s}B_{n+r+k-s}=0, \quad n\in \mathbb{N}, r,s \in \mathbb{Z}^{+} \ \mbox{and $s$ is odd}.
\end{equation}
It can be seen that our approach also applies to the above identity (\ref{kkk}).

Gessel \cite[Theorem 7.3]{gessel03} extended Kaneko's identity
\eqref{kanide} to the following form
\begin{equation}\label{genkan}
\frac{1}{n+1}\sum_{k=0}^{n+1} m^{n+1-k} {n+1\choose k}
\tilde{B}_{n+k}=\sum_{k=1}^{m-1}\bigg((2n+1)k-(n+1)m\bigg)k^n(k-m)^{n-1}.
\end{equation}
Notice  that when $m=1$, the above identity becomes Kaneko's identity.
From the above identity, we obtain the following identity.

\begin{thm} We have
\begin{equation}\label{kanger2}
\frac{1}{(n+3)} \sum_{k=0}^{n+3} m^{n+3-k} {n+3\choose k}
(n+k+3)(n+k+2) \tilde{B}_{n+k}=\sum_{k=1}^{m-1} p(n,m,k)
k^n(k-m)^{n-1},
\end{equation}
where
\begin{align*}
p(n,m,k)=&2(n+2)(2n+3)(2n+5)k^3-2m(n+2)(2n+5)(3n+5)k^2\\[3pt]
&+3m^2(n+2)(2n^2+7n+7)k-m^3(n+1)^2(n+2).
\end{align*}
\end{thm}

\pf Denote the left hand side and the right hand side of
\eqref{genkan} by $L(n,m)$ and $R(n,m)$, respectively. Then we have
\begin{align*}
L(n,m)&=\frac{1}{n+1}\sum_{k=0}^{n+1} m^{n+1-k} {n+1\choose
k}(n+k+1)B_{n+k}\\[3pt]
&=\frac{1}{2\pi i} \oint \frac{1}{e^z-1} \Bigg(\sum_{k=0}^{n+1}
m^{n+1-k} {n+1\choose k}\frac{(n+k+1)}{n+1}
\frac{(n+k)!}{z^{n+k}}\Bigg) dz.
\end{align*}
Denote the summation in the above integral by $S(n,m)$.
By
Zeilberger's algorithm, we find that
\begin{equation}\label{reckan}
z^2S(n+2,m)=2(n+2)(2n+5)S(n+1,m)+m^2(n+1)(n+2)S(n,m).
\end{equation}
Integrating the left hand side of the above recurrence over $z$, we
get \allowdisplaybreaks
\begin{align*}
&\frac{1}{2\pi i} \oint \frac{1}{e^z-1} \big(z^2S(n+2,m)\big) dz\\[3pt]
&=\frac{1}{2\pi i} \oint \frac{1}{e^z-1} \Bigg(\sum_{k=0}^{n+3}
m^{n+3-k} {n+3\choose k}\frac{(n+k+3)}{n+3}
\frac{(n+2+k)!}{z^{n+k}}\Bigg) dz\\[3pt]
&=\frac{1}{(n+3)} \sum_{k=0}^{n+3} m^{n+3-k} {n+3\choose k}
(n+k+3)(n+k+2) \tilde{B}_{n+k}.
\end{align*}
On the other hand, integrating the right hand side of \eqref{reckan}
over $z$ and substituting $L(n,m)$ by $R(n,m)$, we obtain
\begin{align*}
&\frac{1}{2\pi i} \oint \frac{1}{e^z-1}
\bigg(2(n+2)(2n+5)S(n+1,m)+m^2(n+1)(n+2)S(n,m)\bigg) dz\\
&=2(n+2)(2n+5)L(n+1,m)+m^2(n+1)(n+2) L(n,m)\\
&=2(n+2)(2n+5)\sum_{k=1}^{m-1}\bigg((2n+3)k-(n+2)m\bigg)k^{n+1}(k-m)^{n}\\
&\quad+m^2(n+1)(n+2)
\sum_{k=1}^{m-1}\bigg((2n+1)k-(n+1)m\bigg)k^n(k-m)^{n-1}\\
&=\sum_{k=1}^{m-1} p(n,m,k) k^n(k-m)^{n-1},
\end{align*}
as desired. \qed

It is clear that the above identity  \eqref{kanger2} reduces to
Kaneko's identity \eqref{kanger} by setting $m=1$.

\section{Concluding Remarks}

To conclude this paper, we remark that our approach is not
restricted to  identities on Bernoulli and Euler polynomials. It
also applies to sequences $a_0, a_1, a_2, \ldots$  whose  generating
functions $f(z)$ lead to contour integral representations of $a_n$
with hypergeometric integrands. For example, the Genocchi numbers
fall into this framework. We can apply the extended Zeilberger's
algorithm to prove the following identity on Genocchi numbers
\begin{equation}
\sum_{k=0}^{n}( -1) ^{k}{n\choose k}G_{m+k}=\sum_{k=0}^{m}( -1) ^{k}{m\choose k}\sum_{j=0}^{n+k}( -1) ^{j}{n+k\choose j}G_{j},
\end{equation}
where $m,n \in \mathbb{Z}^{+}$. Recall that the Genocchi numbers can
be defined by the generating function
$$
\sum_{n=1}^\infty G_n \frac{z^n}{n!}=\frac{2z}{e^z+1}.
$$

We note that there are other approaches to proving identities
related to special numbers and functions. For example, Paule and
Schneider \cite{paule03} used Karr's summation  algorithm in
difference fields \cite{karr85} and Zeilberger's algorithm  to prove
harmonic number identities and derive new  identities.  Kauers
\cite{kauers07} gave an algorithm which can be applied to verify
many known  identities on Stirling numbers and to discover new
identities. Stan \cite{stan} applied Wegschaider's mathematica
software package {\tt MultiSum} \cite{weg97}  to deal with
identities related to Poisson integrals. Moreover, the package {\tt
MultiSum} can establish multiple index recurrence relations for the
hypergeometric terms with parameters which can also be established
by using the extended Zeilberger's algorithm.

{\noindent \bf Acknowledgments.} We would like to thank Qinghu Hou,
Manuel Kauers,  Peter Paule, Carsten Schneider and  Doron Zeilberger
for their valuable comments. This work was supported by the 973
Project, the PCSIRT Project of the Ministry of Education, the
Ministry of Science and Technology, and the National Science
Foundation of China.

\end{document}